\documentclass{agtart_a}
\pdfoutput=1

\usepackage{pinlabel}

\title{Euclidean Mahler measure and twisted links}

\author[Daniel Silver]{Daniel S Silver}
\givenname{Daniel S}
\surname{Silver}
\address{Department of Mathematics and Statistics\\
University of South Alabama\\\newline
Mobile, AL  36688-0002\\USA}
\email{silver@jaguar1.usouthal.edu}
\urladdr{}

\author{Alexander Stoimenow}
\givenname{Alexander}
\surname{Stoimenow}
\address{Graduate School of Mathematical Sciences\\
University of Tokyo\\\newline
3-8-1, Komaba\\
Tokyo 153-8914\\Japan}
\email{stoimeno@ms.u-tokyo.ac.jp}
\urladdr{}

\author[Susan Williams]{Susan G Williams}
\givenname{Susan G}
\surname{Williams}
\address{Department of Mathematics and Statistics\\
University of South Alabama\\\newline
Mobile, AL  36688-0002\\USA}
\email{swilliam@jaguar1.usouthal.edu}
\urladdr{}

\volumenumber{6}
\issuenumber{}
\publicationyear{2006}
\papernumber{21}
\startpage{581}
\endpage{602}

\doi{}
\MR{}
\Zbl{}


\received{26 March 2005}
\revised{}
\accepted{15 March 2006}
\published{7 April 2006}
\publishedonline{7 April 2006}
\proposed{}
\seconded{}
\corresponding{}
\editor{}
\version{}

\arxivreference{math.GT/0412513}

\AtBeginDocument{\let\bar\wbar\let\tilde\wtilde\let\hat\what}
\numberwithin{equation}{section}

\makeatletter
\def\cnewtheorem#1[#2]#3{\newtheorem{#1}{#3}[section]
\expandafter\let\csname c@#1\endcsname\c@thm}

\newtheorem{thm}{Theorem}[section]   
\cnewtheorem{lem}[thm]{Lemma}         
\cnewtheorem{cor}[thm]{Corollary}    
\cnewtheorem{prop}[thm]{Proposition}  

\theoremstyle{definition}
\cnewtheorem{defn}[thm]{Definition}   
\cnewtheorem{rem}[thm]{Remark}        
\cnewtheorem{exa}[thm]{Example}   

\makeautorefname{defn}{Definition}

\makeatother  
\def\a{{\alpha}}

\def\e{{\epsilon}}

\def\d{{\delta}}

\def\Z{{\mathbb Z}}

\def\<{{\langle}}
\def\>{{\rangle}}

\def\D{{\Delta}}

\def\i{{\iota}}
\def\R{{\mathbb R}}
\def\<{{\langle}}
\def\>{{\rangle}}

\def\m1{{\quad ({\rm mod\ 1})}}

\def\={\ {\buildrel \cedot \over =}\ }
\def\tsp{{trivial split component}}
\def\tsps{trivial split components}

\begin{document}

\begin{asciiabstract}
 If the twist numbers of a collection of oriented alternating link
 diagrams are bounded, then the Alexander polynomials of the
 corresponding links have bounded euclidean Mahler measure (see
 Definition 1.2).  The converse assertion does not hold. Similarly, if
 a collection of oriented link diagrams, not necessarily alternating,
 have bounded twist numbers, then both the Jones polynomials and a
 parametrization of the 2-variable Homflypt polynomials of the
 corresponding links have bounded Mahler measure.
\end{asciiabstract}

\begin{webabstract}
 If the twist numbers of a collection of oriented alternating link
 diagrams are bounded, then the Alexander polynomials of the
 corresponding links have bounded euclidean Mahler measure (see
 Definition 1.2).  The converse assertion does not hold. Similarly, if
 a collection of oriented link diagrams, not necessarily alternating,
 have bounded twist numbers, then both the Jones polynomials and a
 parametrization of the 2--variable Homflypt polynomials of the
 corresponding links have bounded Mahler measure.
\end{webabstract}

\begin{abstract}  
 If the twist numbers of a collection of oriented alternating link
 diagrams are bounded, then the Alexander polynomials of the
 corresponding links have bounded euclidean Mahler measure (see
 \fullref{emm}).  The converse assertion does not hold. Similarly, if
 a collection of oriented link diagrams, not necessarily alternating,
 have bounded twist numbers, then both the Jones polynomials and a
 parametrization of the 2--variable Homflypt polynomials of the
 corresponding links have bounded Mahler measure.
\end{abstract}

\maketitle

\section{Introduction} 
\label{intro}
If $f(t)$ is a nonzero polynomial with complex coefficients,
$$f(t) = b\cdot \prod_{i=1}^n (t-\a_i) \in \mathbb C[t],$$
then its {\it Mahler measure} \cite{mahler} is
$$M(f) = |b|\cdot\prod_{i=1}^n \max\{|\a_i|, 1\}.$$
By convention, the Mahler measure of the zero polynomial is defined to
be $0$. A well-known theorem of Kronecker implies that a monic (ie,
$|b|=1$) integer
polynomial has Mahler measure 1 if and only if it is a product of monomials and cyclotomic polynomials.

In 1933, D\,H Lehmer discovered a monic integer polynomial of degree 10, 
$$L(t) = t^{10}+t^9 -t^7-t^6-t^5-t^4-t^3+t+1,$$
with a single zero of modulus greater than 1, equal approximately to 1.17628 \cite{lehmer}. He asked whether given any $\epsilon >0$ there exists an integer polynomial $f(t)$ such that
$1 < M(f)< 1+\epsilon$. Lehmer's question remains open. In fact, no integer polynomial has been found with Mahler measure less than $M(L)$ but greater than 1. 

Jensen's formula motivates a definition of Mahler measure for polynomials of more than one variable \cite{mahler}. For any nonzero polynomial $f(t_1, \ldots, t_d)$ with complex coefficients, the Mahler measure of $f$ is defined by the following integral, which is possibly singular but nevertheless convergent:
$$M(f)=\exp \int_0^1\cdots \int_0^1 \log |f(e^{2\pi i \theta_1}, \ldots, e^{2\pi i \theta_d})|d\theta_1\cdots d\theta_d.$$
Lehmer's question for polynomials of any higher degree is equivalent
to the question for 1--variable polynomials (Boyd \cite{boyd}). The reader who
is interested in additional background information on Mahler measure is
encouraged to consult Everest and Ward \cite{EW}.

\begin{rem}\label{rem1.1}
The Mahler measure of a nonzero Laurent polynomial 
$t^{-r}f(t)$, $f(t)\in {\mathbb C}[t]$, is defined to be $M(f)$. Likewise, we define the Mahler measure of $t^{1/2}f(t)$ to be $M(f)$. Thus  the Mahler measure of both the 1--variable Alexander and the Jones polynomial of any oriented knot or link is well defined. 
Similarly one deals with several variable Laurent polynomials. \end{rem}

A topological interpretation of Mahler measures of single- and
multivariable-Alexander polynomials of links in terms of homology
growth of branched covers was given by Silver and Williams \cite{SW02}. Moreover,
the Mahler measure of the Alexander polynomial of a fibered link
is a suitably defined growth rate of Lefschetz numbers of the fiber
monodromy \cite{SW04a}. Slowing down the
dynamics of pseudo-Anosov homeomorphisms is equivalent to Lehmer's
question. In this sense, it is a question about fibered knots and
links and the dynamical properties of their monodromies.

Lehmer's polynomial $L(t)$ occurs (up to the interchange $t
\leftrightarrow -t$) as the Alexander polynomial of a knot. Perhaps
the most interesting such knot is the $(-2,3,7)$--pretzel knot, a
fibered hyperbolic knot with noteworthy properties (see for example
Gordon \cite{gordon}, Hironaka \cite{hiro}). If Lehmer's question has
an affirmative answer, then the polynomial $f(t)$ can be chosen to be
a (1--variable) Alexander polynomial of a fibered hyperbolic
arborescent (or Conway-algebraic) 2--component link in $S^3$
(Stoimenow \cite{stoi05}) or a fibered hyperbolic knot in a lens space
$L(n,1)$ \cite{SW04a}. (For a knot in $S^3$ one needs one of the
conditions $f(1)= \pm 1$, whose meaning in Lehmer's question is
unclear so far.)

Callahan, Dean and Weeks argue in  \cite{CDW} that for hyperbolic knots, the adjective ``simple" should reflect geometric properties. Consequently, they propose that a hyperbolic knot or link should be considered simple if its complement can be constructed with relatively few ideal tetrahedra. Examples in \cite{SW04} suggest that the multivariable Alexander polynomial of such a knot or link has small Mahler measure.  Examples of simple hyperbolic knots in  \cite{CKP} suggest a similar statement for the Jones polynomial $V_\ell$. The $(-2,3,7)$--pretzel knot is one example, with complement composed of only 3 ideal tetrahedra. 

Often the logarithm of Mahler measure arises as the topological entropy of algebraic dynamical systems. In \cite{LW} Lind and Ward defined $p$--adic and euclidean entropy for automorphisms of generalized solenoids, the type of algebraic dynamical systems that arise in \cite{SW02}, and they showed that topological entropy is their sum. The two components correspond to the two contributions  in the definition of Mahler measure of a polynomial $f(t)$, one from the leading coefficient of $f$ and the other from its zeros. This motivates the following.

\begin{defn}\label{emm}
If $f(t)\in {\mathbb C}[t]$ has zeros $\a_1, \ldots, \a_n$, then its
{\it euclidean Mahler measure} is
$$M_e(f) = \prod_{i=1}^n \max\{|\a_i|, 1\}.$$
\end{defn}

\begin{rem}\label{rem1.3} 
Euclidean Mahler measure of the Alexander polynomial of a knot $k$ is
a natural quantity.  The homology group $H_1(\tilde X; {\mathbb R})$
of the infinite cyclic cover of $k$ with real coefficients is a
finite-dimensional vector space, and a generator of the deck
transformation group of $\tilde X$ induces an automorphism. It is not
difficult to see that the product of moduli of those eigenvalues that
are outside the unit circle coincides with $M_e(\D_k)$. \end{rem}

Our main result is concerned with sequences of knots or links that
have 1--variable Alexander polynomials with euclidean Mahler measures
tending to infinity or Jones polynomials with Mahler measures tending
toward infinity.  In such cases we can conclude that the twist numbers
also tend toward infinity. In the case that the knots or links are
alternating and hyperbolic, their volumes increase without
bound. (Dasbach and Lin \cite{DL} give an expression of the twist
number of alternating diagrams in terms of the Jones polynomial, which
gives a different relation between twist number and volume.)

\medskip
{\bf Acknowledgements}\qua The authors are grateful to Hiroshima
University and Osaka City University for their hospitality while this
work was performed. Their visits were funded by the program
``Constitution of wide-angle mathematical basis focused on knots,"
directed by Professor Akio Kawauchi and part of the 21st Century COE
Program.  The calculations were aided by Knot, a program developed by
K Kodama \cite{koda}.  The first and third authors were partially
supported by NSF grant DMS-0304971, and the second author was
supported by JSPS Postdoc grant P04300.


\section{Statement of main theorem} 
By {\it full-twisting} an oriented link $\ell$, we mean cutting a pair
of adjacent arcs of a diagram $D$ for $\ell$, inserting some number
$q$ of full twists (right-handed if $q$ is positive and left-handed
otherwise) and then reattaching the arcs. In this way, we obtain a
sequence of links $\ell_q$. As $q$ goes to infinity, the Mahler
measures of the Jones polynomials $V_{\ell_q}$ converge (Champanerkar
and Kofman \cite{CK}); the Mahler measures of the Alexander polynomials
$\D_{\ell_q}$ will also converge, provided that the twisted arcs are
coherently oriented \cite{SW04}. If we insert arbitrary numbers $q_i$
of full-twists at several sites, then as all $q_i$ grow without bound,
we again we have convergence. However, the set of Mahler measures
produced will generally have infinitely many distinct limit points
(cf \fullref{ex5.1}). \fullref{main} implies that the limit points are
bounded.

In the case that the arcs being twisted are not coherently oriented,
the Mahler measures of the Alexander polynomials of $\ell_q$ can grow
without bound. Twist knots provide simple examples.  However, we will
see that if the diagram $D$ is alternating and we replace Mahler
measure by euclidean Mahler measure, then again limits exist and limit
points are bounded.

In order to state the main result, we need the following notions. 

Let $D$ be a link diagram. A {\it bigon region} is a complementary region with exactly two crossings in its boundary. A {\it twist} is either a connected row of bigon regions, maximal in the sense that it is not part of a longer row of bigons, or else it is a single crossing that is adjacent to no bigon region. The number of twists of $D$ is called the {\it twist number}, denoted by $t(D)$. The twist number of a link $\ell$ is the  minimal twist number of all of its diagrams. 

A {\it \tsp{}} of a link diagram 
is a circle (with no self-crossings) disjoint from the rest of the 
diagram. A {\it \tsp{}} of a link is an unknotted component separable 
by a hyperplane from the rest of the link.

We define the {\it length} $||p||$ of a multivariable polynomial $p$
to be the sum of the absolute values of its coefficients. It follows
from the triangle inequality and monotonicity of the log function that
$M(p)\leq ||p||$. (The inequality for single-variable polynomials
appears in Everest and Ward \cite{EW} as an exercise. For another
argument, see \cite[Lemma 6.1]{stoi03} for $p\in \Z[x]$. The same
proof works for $p\in \R[x]$, and as remarked there also for
polynomials in several variables.)

Let $P_\ell(v,z)$ be the {\it Homflypt skein polynomial} with the
skein relation
$$v^{-1}P_{\ell_+}(v,z)-vP_{\ell_-}(v,z) = zP_{\ell_0}(v,z),$$
normalized in the usual way so that the polynomial of the trivial
knot is $1$. We will consider a $2$--variable parametrization
$$\hat P _\ell(v,t)\,=\,P_\ell(v, t^{1/2}-t^{-1/2})$$
that is related by a variable change
preserving Mahler measure and length to the polynomial
$X(q,\lambda)$ that appears in Jones \cite{jones}. The {\it Jones polynomial}
is obtained by setting $q=\lambda=t$, or in other words:
$$V_\ell(t)\,=\,\hat P_\ell(t,t)\,=\,P_\ell(t,t^{1/2}-t^{-1/2})\,,$$
and the {\it Alexander polynomial} (with a particular choice of
normalization) by
$$\D_\ell(t)\,=\,\hat P_\ell(1,t)\,=\,P_\ell(1,t^{1/2}-t^{-1/2})\,.$$
Here and below the term ``Alexander polynomial" will refer to the {\it
one}-variable version; its several-variable relative will be referred
to as the ``multi-variable Alexander polynomial." In each case, the polynomial 
is defined up to multiplication by a unit. 

\begin{thm}\label{main}
{\rm(1)}\qua If ${\cal D}$ is a set of oriented link
diagrams with twist numbers at most $n$, and no \tsps, then the
sets of lengths
$$\{||(t+1)^nV_{\ell(D)}(t)||\mid D \in {\cal D}\}$$
and
$$\{\,||(v^2-1)^n(t+1)^n\hat P_{\ell(D)}(v,t)||\mid D \in 
{\cal D}\}$$ are bounded.

{\rm(2)}\qua If ${\cal D}$ is a set of oriented alternating link diagrams with
bounded twist numbers, then the set of euclidean Mahler measures
$\{M_e(\D_{\ell(D)}) \mid D \in {\cal D}\}$ is  bounded. \end{thm}

\begin{cor}\label{jonescor}
If ${\cal D}$ is a set of oriented link diagrams with bounded twist
number, the set of Mahler measures $\{M(V_{\ell(D)}) \mid D \in {\cal
D}\}$ and $\{M(\hat P_{\ell(D)}) \mid D \in {\cal D}\}$ of the Jones
and parametrized Homflypt polynomials are bounded. \end{cor}

\begin{proof} \fullref{jonescor} follows from the fact that the Mahler measure of a polynomial is bounded by its length, while multiplying a polynomial by a cyclotomic polynomial does not change its Mahler measure. \end{proof}

\begin{cor}\label{cor2.3} If $\ell_q$ is a sequence of prime alternating 
hyperbolic links such that either $M(V_{\ell_q})$ or $M_e(\D_{\ell_q})$ increases without bound, then ${\rm Volume}(S^3\setminus \ell_q)$ also increases without bound.\end{cor}

\begin{proof} Each link $\ell_q$ admits a connected alternating
diagram $D_q$ that is prime in the sense that any simple closed curve in the
plane meeting $D_q$ in exactly two points disjoint from the crossings bounds
a region containing no crossings. Theorem 1 of Lackenby \cite{lacken} implies that
$$v_3 (t(D_q)-2)/2 \le {\rm Volume}(S^3\setminus \ell_q)\,,$$
where $v_3\ (\approx  1.01494)$ is the volume of a regular hyperbolic ideal 3--simplex. 
Since the twist numbers $t(D_q)$ increase without bound, by \fullref{main},  the volumes of $S^3\setminus \ell_q$ also increase without bound. \end{proof}


\section{Twisting} We prepare for the proof of \fullref{main} by establishing results about the Alexander, Jones and skein polynomials of an
oriented link when twists at several sites are performed.

\begin{figure}[ht!]\small
\labellist
\pinlabel (1) at  78 35
\pinlabel (2) at 244 35
\pinlabel (3) at 403 35
\endlabellist
\begin{center}
{\includegraphics[width=3in]{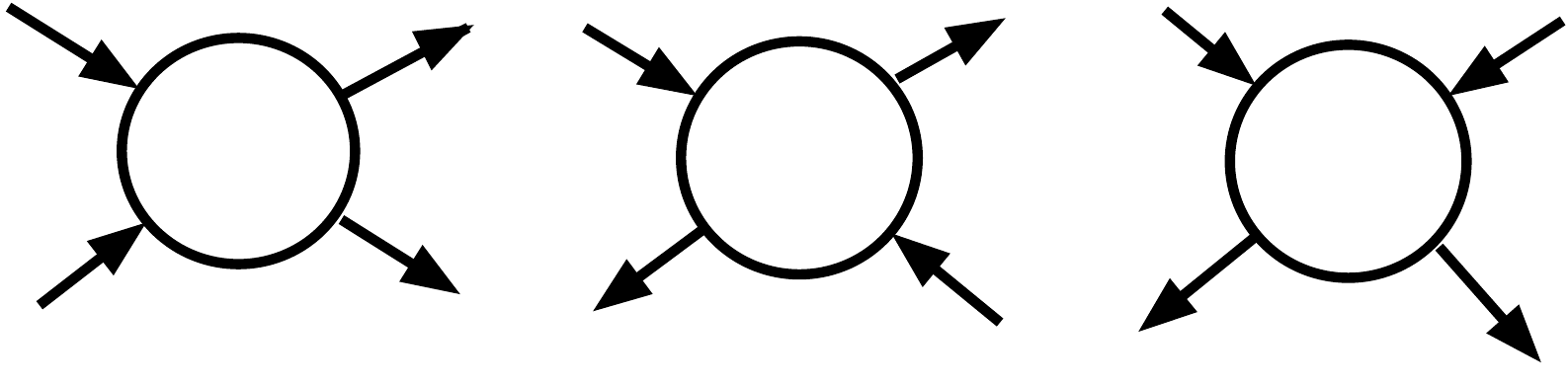}} 
\caption{Vertex types of a wiring diagram} 
\label{twist1}
\end{center}
\end{figure}

By a {\it wiring diagram} we mean a planar diagram consisting of $n$
vertices or {\it twist sites}  $v_1, \ldots, v_n$ of the three possible types shown in \fullref{twist1}, connected by oriented arcs. A vertex of the first type is a {\it parallel twist site} while a vertex of the second or third type is {\it anti-parallel}.
We call the number of twist sites the {\it order} of the wiring diagram. 

Given a wiring diagram of order $n$, we consider $n$--tuples $(q_1, \ldots, q_n)$, where $q_i \ne \infty$ if $v_i$ is of the first type; $q_i$ is even if $v_i$ is of the second type; and $q_i$ is odd if $v_i$ is of the third type. We obtain an oriented link $\ell(q_1, \ldots, q_n)$ by replacing the vertex $v_i$ by $T_{-q_i}$ if $v_i$ is of the first type, and $T_{q_i}$ if it is of the second or third type. (This choice assures that the sign of $q_i$ agrees with the sign of the crossings at that site in the resulting link.) We will use the previous
term {\it twist} to refer to each embedded tangle $T_{q_i}$. Clearly any oriented link diagram with twist number $n$ can be described by a wiring diagram of order $n$. Note that a twist consisting of a single crossing can be regarded as either parallel or anti-parallel, depending on the choice of wiring diagram.

\begin{figure}[ht!]\small
\labellist
\pinlabel {$q<0$} at 146 239
\pinlabel {$q=0$} at 429 239
\pinlabel {$q>0$} at 146 30
\pinlabel {$q=\infty$} at 429 30
\endlabellist
\begin{center}
{\includegraphics[width=2in]{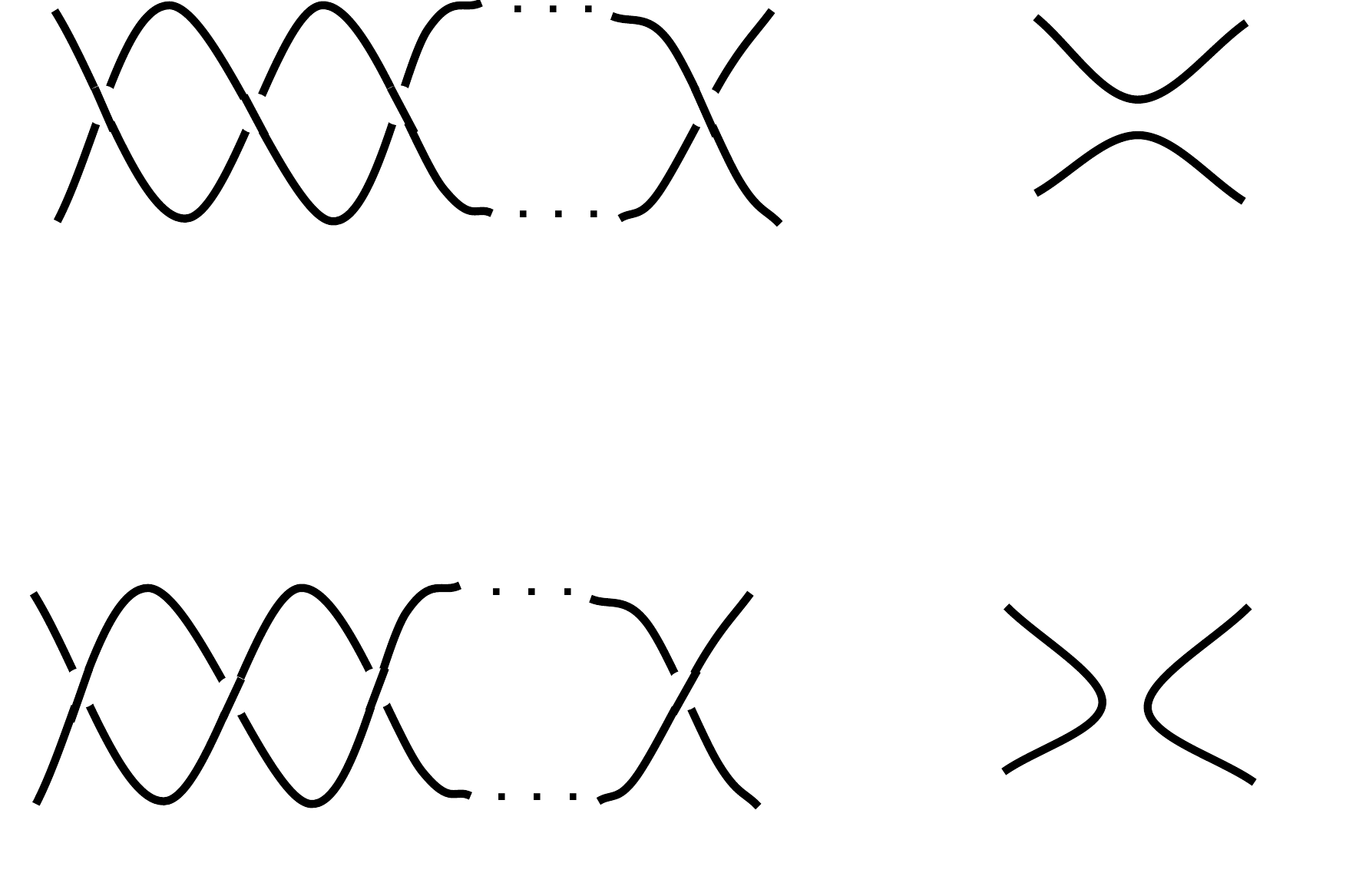}} 
\caption{Tangles $T_q$} 
\label{twist2}
\end{center}
\end{figure}
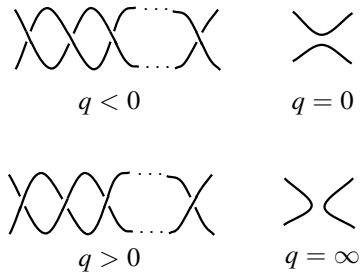

\begin{lem}\label{lem3.1} 
Let $D$ be an oriented link diagram, and let $D_q$ (for $q\in {\mathbb
Z}$) be the diagram that results from replacing a pair of adjacent
parallel arcs of $D$ by the tangle $T_{-q}$ in such a way that
$D_0=D$. Let $\ell_q$ be the links described by $D_q$. Then
$$(t+1)\D_{\ell_q}(t) = \Big(\D_{\ell_0}(t) + t^{1/2}\D_{\ell_1}(t)\Big)
t^{q/2} + (-1)^q \Big(t\D_{\ell_0}(t) - t^{1/2}\D_{\ell_1}(t)\Big)
t^{-q/2},$$
$$(t+1)V_{\ell_q}(t) = \Big(V_{\ell_0}(t) + t^{-1/2}V_{\ell_1}(t)\Big)t^{3q/2} +(-1)^q \Big(tV_{\ell_0}(t)- t^{-1/2}V_{\ell_1}(t)\Big)t^{q/2},$$
$$(t+1) \hat P_{\ell_q}(v,t)= v^q\Big[\Big(\hat P_{\ell_0}
(v,t)+ v^{-1}t^{1/2}\hat P_{\ell_1}(v,t)\Big)t^{q/2} +$$
$$ (-1)^q\Big(t\hat P_{\ell_0}(v,t) - v^{-1}t^{1/2}\hat P_{\ell_1}(v,t)\Big)t^{-q/2}\Big].$$
\end{lem}

\begin{proof} Let $D_q$ be the diagram with $q$ twists. Set
$$p_q(v,t) = v^{-\omega(D_q)}\hat P_{D_q}(v, t),$$
where
$\omega(D_q)$ is the writhe. From the skein relation
$$v^{-1}P_+(v,z)-vP_-(v,z) = zP_0(v,z),$$
we have 
\begin{equation}p_{q+2}(v,t)-p_q(v,t)=(t^{1/2}-t^{-1/2})p_{q+1}(v,t).
\label{eq3.1}
\end{equation}
For the case $q\ge 0$ we use the generating function
$$f(v,t,x)=\sum_{q=0}^\infty p_q(v,t)x^q.$$
Multiplying equation \ref{eq3.1} by $x^q$ and summing from $0$ to $\infty$ gives
$${{f(v,t,x)-p_0(v,t)-p_1(v,t)x}\over x^2} - f(v,t,x) = (t^{1/2}-t^{-1/2}){{f(v,t,x)-p_0(v,t)}\over x},$$
which can be written as
$$f(v,t,x)= {{p_0(v,t)+\Big(p_1(v,t)-(t^{1/2}-t^{-1/2})p_0(v,t)\Big)x}\over{1-(t^{1/2}-t^{-1/2})x-x^2}}.$$
A partial fraction decomposition gives
$$f(v,t,x)={1\over t+1}\Big[ {{p_0(v,t)+t^{1/2}p_1(v,t)}\over{1-t^{1/2}x}} + {{tp_0(v,t)-t^{1/2}p_1(v,t)}\over{1+t^{-1/2}x}}\Big].$$
Taking series expansions in $x$ of these rational functions, and equating coefficients, shows that 
$$p_q(v,t)= {\textstyle{1\over t+1}}\Big[\Big(p_0(v,t)+t^{1/2}p_1(v,t)\Big)t^{q/2} + (-1)^q\Big(tp_0(v,t)-t^{1/2}p_1(v,t)\Big)t^{-q/2}\Big].$$
Since $\omega(D_q) = q+\omega(D_0),$ we see the statement for $\hat P$.
The remaining results follow by setting $v=1$ for the
Alexander polynomial, and setting $v=t$ for the Jones
polynomial. The case $q<0$ is handled in the same manner, using the generating  function $g(v,t,x)= \sum_{r=0}^\infty p_{1-r}(v,t)x^r$ and then setting $r=1-q$. \end{proof}

\begin{rem} (1)\qua Champanerkar and Kofman showed in \cite{CK} using Jones-Wenzl idempotents that for full twists, $V_{\ell_q}(t)$ can be expressed as a rational function of $t$ and $t^q$. In fact, their result holds more generally for $q$ full twists on any number of strands (with arbitrary orientation).  

(2)\qua A special case of the first formula of \fullref{lem3.1} appears in Bhatty \cite{bhatty}. \end{rem}

\begin{lem}\label{lem3.3} Let $D$ be an oriented link diagram, and let $D_q$ (for $q\in 2{\mathbb Z}\cup\{\infty\}$) be the diagram that results from replacing a pair of adjacent anti-parallel arcs of $D$  by the tangle $T_q$ in such a way that $D_0=D$.
Let $\ell_q$ be the link described by $D_q$. Then for $q \in 2{\mathbb Z}$, 
$$\D_{\ell_q}(t) = \D_{\ell_0}(t)+ {q\over 2}(t^{1/2}-t^{-1/2})\D_{\ell_\infty}(t), $$
$$(t+1)V_{\ell_q}(t)= \biggr((t+1)V_{\ell_0}(t)+ t^{1/
2}V_{\ell_\infty}(t)\biggr)t^q - t^{1/2}V_{\ell_\infty}(t),$$
$$ P_{\ell_q}(v,z)\,=\,{{v^{q}-1}\over{v-v^{-1}}}z\,P_{\ell_\infty}(v,z)+v^qP_{\ell_0}(v,z)\,.$$
\end{lem}

\fullref{lem3.3} is a straightforward consequence of the Homflypt skein
relation. (See formula (8) in \cite{stoi98}, or formula (4) in 
\cite{stoi00}, including a correction of the misprint in
the first reference.) Note that the case of an odd number of
anti-parallel twists can be handled by letting $D=D_0$ have
a single crossing at the twist site in question.

\begin{thm}\label{th3.4} 
Consider a wiring diagram of order $n$, and let $\ell(q_1, \ldots,
q_n)$ be the oriented link defined at the beginning of the
section. Assume that sites $1,\ldots, m$ are anti-parallel, and the
remaining sites are parallel.

{\rm(1)}\qua The Jones polynomial of $\ell(q_1, \ldots, q_n)$ satisfies:
$$(t+1)^{n}V_{\ell(q_1, \ldots, q_n)}(t) =
W(t^{1/2}, t^{q_1/2}, \ldots, t^{q_n/2}),$$
where $W(u, w_1, \ldots, w_n) \in {\mathbb Z}[u^{\pm1}, w_1^{\pm1},
\ldots, w_n^{\pm1}]$ depends only on the wiring diagram describing
$\ell(q_1, \ldots, q_n)$ and the parities of $q_1, \ldots, q_n$.

{\rm(2)}\qua The Alexander polynomial of $\ell(q_1, \ldots, q_n)$ satisfies: 
$$(t+1)^{n-m}\D_{\ell(q_1, \ldots, q_n)}(t) = \sum {{q_1^{\d_1}\cdots
q_m^{\d_m}}\over 2^m} X_{\d_1, \ldots, \d_n}(t^{1/2}, t^{q_{m+1}/2},
\ldots, t^{q_n/2}),$$
where the sum ranges over all choices of $\d_i \in \{0,1\}$, $i=1,
\ldots, n$, and each $X_{\d_1, \ldots, \d_n}\in {\mathbb Z}[u^{\pm1},
w_1^{\pm1}, \ldots, w_{n-m}^{\pm1}]$ depends only on the wiring
diagram describing $\ell(q_1, \ldots, q_n)$ and the parities
of the $q_1, \ldots q_n$.

{\rm(3)}\qua The Homflypt polynomial of $\ell(q_1, \ldots, q_n)$
satisfies: 
\begin{equation*}\begin{split}
(t+1)^{n-m}&(v^2-1)^{m}\hat P_{\ell(q_1, \ldots, q_n)}(v,t)\  \\
&=\ Y(v^{1/2}, t^{1/2}, v^{q_1/2}, \ldots, v^{q_m/2},
t^{q_{m+1}/2}, \ldots, t^{q_n/2}),\end{split}\end{equation*}
where $Y(u,w,u_1,\ldots,u_m,w_{m+1}, \ldots, w_n) \in {\mathbb Z}[
u^{\pm1},w^{\pm1},\ldots, w_n^{\pm1}]$ depends only on the wiring
diagram describing 
$\ell(q_1, \ldots, q_n)$ and the parities of $q_1, \ldots, q_n$.\end{thm}

\begin{proof} 
Using either \fullref{lem3.1} (if $v_n$ is parallel) or
\fullref{lem3.3} (if $v_n$ is anti-parallel), we can write
$V_{\ell(q_1, \ldots, q_n)}$ as a polynomial in $t^{1/2}, t^{q_n/2}$
and the expressions $V_{\ell(q_1, \ldots, q_{n-1},0)}(t), V_{\ell(q_1,
\ldots, q_{n-1},1)}(t)$ and $V_{\ell(q_1, \ldots, q_{n-1},\infty)}(t)$
(only two of the expressions will appear). The polynomial depends only
on the parity of $q_n$. Continuing in this fashion with each twist
site in turn, we write $V_{\ell(q_1, \ldots, q_n)}$ as a polynomial in
$t^{1/2}, t^{q_1/2}, \ldots, t^{q_n/2}$ and the polynomials
$V_{\ell(\e_1, \ldots, \e_n)}(t) \in {\mathbb Z}[t^{\pm 1/2}]$, where
$\e_i \in \{0,1,\infty\}.$ This polynomial has the desired form, and
depends only on the wiring diagram and the parities of $q_1, \ldots,
q_n$.  The argument for the Homflypt polynomial is similar.

We prove the second assertion of \fullref{th3.4}. For $i=1, \ldots, m$, we set $\d_i=1$ if $\e_i=\infty$ and $\d_i=0$ if $\e_i=0$ ($q_i$ even) or $\e_i=1$ ($q_i$ odd). Using Lemma 3.3 we can write 
$$\D_{\ell(q_1, \ldots, q_n)}(t)= \sum \bigg \lfloor {q_1\over 2}\bigg \rfloor^{\d_1}\cdots \bigg \lfloor {q_m\over 2}\bigg \rfloor^{\d_m}\D_{\ell(\e_1, \ldots, \e_m, q_{m+1}, \ldots, q_n)}(t).$$
Now we use \fullref{lem3.1} as before to write $\D_{\ell(\e_1, \ldots, \e_m, q_{m+1}, \ldots, q_n)}(t)$ as a linear combination of polynomials
in $t^{1/2}, t^{q_{m+1}/2}, \ldots, t^{q_n/2}$, that depend on the
parity of $q_{m+1},\ldots,q_n$. \end{proof}

\fullref{th3.4} has independent interest. It says that if a collection
of links has bounded twist number (and no \tsps), then the result of
multiplying either the Jones or Alexander polynomials by a fixed
polynomial (a power of $t+1$) is a collection of polynomials with a
bounded number of nonzero coefficients. \fullref{lem3.5} shows that a
collection of polynomials $fg^N$ will have such a form only in trivial
cases.

\begin{lem}\label{lem3.5} If $f(t), g(t)\in {\mathbb Z}[t^{\pm 1}]$, $f(t)\ne 0$, $g(t)$ not a unit, then
the number of nonzero coefficients of $fg^N$ tends to infinity as $N$ increases without bound.\end{lem}

\begin{proof} We may assume without loss of generality that $f(t)= 1+a_1t+\cdots$ and $g(t)= 1 +b_mt^m+\cdots,$ where $a_i, b_i \in {\mathbb Z}$ and $b_m\ne 0$. We show that for each $k\ne 0$, the coefficient $c_{km}^{(N)}$ of $t^{km}$ in $fg^N$ is nonzero for all sufficiently large $N$. This coefficient is  given by the  sum
$$c_{km}^{(N)} = \sum_{n+n_1+\cdots + n_N=km} a_n  b_{n_1}\cdots b_{n_N},$$
with $n,n_i\ge 0$.
Splitting the sum over the number $p$ of nonzero $n_i$, we can write it as
$$c_{km}^{(N)} = \sum_{p=1}^k {N\choose p}\sum a_n  b_{n_1}\cdots b_{n_p},$$
where the second summation is taken over all $n, n_1, \ldots, n_p$ with 
$m\le n_i, i=1, \ldots, p$, and $n+n_1+\cdots +n_p=km$. The second summation  is independent of $N$ for each $p$, so for $N\ge k$ we may regard $c_{km}^{(N)}$ as a polynomial in $N$. The leading term of the polynomial comes from the summand corresponding to $p=k$, and it is given by ${1\over k!}b_m^kN^k$. Hence $c_{km}^{(N)}$ is nonzero for $N$ sufficiently large. \end{proof}

\begin{cor} If $\ell$ is a link with non-trivial Jones
or Alexander polynomial (ie, not equal to $1$) and no \tsps, then 
the twist number of the connected sum $\ell\ \sharp\  \ell\ \sharp\  \cdots\  \sharp\  \ell$ tends to infinity as the number of summands increases without bound.\end{cor}

\begin{proof}Let $g$ denote either the Jones or Alexander polynomial of $\ell$. Since $g$ is non-trivial, it is not a unit. For the Alexander polynomial, this is a consequence of normalization, as a specialization of $P$ (see paragraph preceding \fullref{main}). In the case of the Jones polynomial, it follows from the properties given in section 12 of \cite{jones}.

The Jones and Alexander polynomials of $\sharp_{i=1}^N \ell$ have the form $g^N$. If the connected sums have bounded twist numbers, then after multiplying by a suitable power $f(t)$ of $t+1$, the polynomials $fg^N$ have a bounded number of nonzero coefficients, contradicting \fullref{lem3.5}. \end{proof}

The following proposition will be needed for the proof of \fullref{main}.

\begin{prop} Let $f(x, y_1, \ldots, y_n) \in {\mathbb Z}[x^{\pm1}, y_1^{\pm1}, \ldots, y_n^{\pm1}]$. There is a constant $C$ with 
$$||f(x, x^{q_1}, \ldots, x^{q_n})|| \le C,$$
for all $(q_1, \ldots, q_n)\in {\mathbb Z}^n.$\end{prop}

\begin{proof} Write $f(x, y_1, \ldots, y_n)$ as the sum of terms
$f_j(x)g_j(y_1, \ldots, y_n)$, $1\le j\le m$, where $g_j$ is a product of powers of the $y_i$. Then  
$$||f_j(x)g_j(x^{q_1}, \ldots, x^{q_n})|| = ||f_j(x)||.$$
Hence we have
$$||f(x, x^{q_1}, \ldots, x^{q_n})|| \le \sum_{j=1}^m ||f_j(x)||,$$
from the triangle inequality. \end{proof}

The proof of \fullref{main} also requires some facts about Murasugi 
products.
Consider an oriented diagram $D$ for a link $\ell$. 
By smoothing crossings according to the well-known algorithm of Seifert
\cite{seifert}, we obtain a number of (possibly nested) {\it Seifert
disks} in the plane.  An orientable spanning surface for the link,
called a {\it canonical Seifert surface},
is then seen by connecting the disks with half-twisted
bands corresponding to the crossings of the diagram. 

The boundary of a Seifert disk is called a {\it Seifert circle}.
If $C$ is a Seifert circle, then it decomposes the plane into two
closed regions $U, V$ meeting along $C$. We say that $C$ is {\it
separating} if both $(U\setminus C)\cap D$ and $(V\setminus C)
\cap D$ are non-empty, otherwise $C$ is {\it non-separating}.
If $C$ is separating, then let $D_1$ and $D_2$ be the diagrams
constructed from $D\cap U$ and $D\cap V$, filling gaps with arcs
from $C$ where they are needed. We say that $D$ is a *-{\it product}
(or {\it Murasugi product}) of $D_1$ and $D_2$, and write $D= D_1*D_2$. 

A diagram $D$ is {\it special} if it does not decompose as a Murasugi
product, in other words, if and only if it has no separating Seifert
circle. A general oriented diagram $D$ can be decomposed along its
separating Seifert circles into a product
$D_1*\cdots *D_r$ of special diagrams.
We call $D_i$ the {\it special (Murasugi) factors} of $D$.
If $D$ is alternating, so
are the $D_i$; they are {\it special alternating diagrams}.
Notice that each twist of $D$ is contained in some factor. 
The diagram $D$ is connected if and only if each $D_i$ is
connected. For additional background,
the reader is advised to consult Murasugi \cite{mura60,mura62} or
Cromwell \cite{cromwell}.

We shall assume in the following that $D$ is a connected
oriented alternating diagram, and $D_1*\cdots *D_r$ its decomposition
into special alternating diagrams $D_i$. Let 
$\ell, \ell_1, \ldots, \ell_r$ be the links represented by
$D, D_1,\ldots, D_r$, respectively. We will make use of the fact, noted in \cite{mura60}, that the leading coefficient of $\D_\ell(t)$, up to sign, is the product of the leading coefficients of $\D_{\ell_1}(t), \ldots, \D_{\ell_r}(t)$. The fact holds more generally when $D$ is a homogeneous diagram (Murasugi and Przytycki \cite{MP}). 

Associated to an oriented alternating diagram $D$ for a link $\ell$ there is a
graph $\Gamma$ obtained in the following manner.  Checker-board color the
regions of the diagram with black and white, and let the vertices of
$\Gamma$ correspond to the black regions. Two vertices are connected
by an edge for each crossing shared by the corresponding regions.
Note that $\Gamma$ is planar, (generally) with multiple edges,
and interchanging colors transforms $\Gamma$ into the dual graph.

When $D$ is a connected special alternating diagram, it is possible to
checker-board color so that each white region has a Seifert circle as
its boundary. Then the degree of each vertex of $\Gamma$ is even, and we can orient the edges in such a way that they alternate in and out as we travel around each vertex. Such an orientation is unique up to global orientation reversal. (Note that the graph is connected since by assumption $D$ is connected.) 
Fix a vertex $v$, called the {\it root}. The edges of any spanning tree $T\subset \Gamma$ can be uniquely oriented ``toward the root" so that each vertex other than $v$ has exactly one outgoing edge. 

Define $\i(T)$ to be the number of edges of $T$ that disagree with the
orientation of $\Gamma$. We call such edges {\it incoherent}, and the
others {\it coherent}. We call $T$ incoherent (resp.\ coherent) if all
its edges are incoherent (resp.\ coherent).  In Murasugi and Stoimenow
\cite{MS} it is shown that, up to units in $\Z[t^{\pm 1/2}]$,
\begin{equation}
\D_\ell(-t) = \sum_{T\subset \Gamma} t^{\i(T)},
\label{eq3.2}
\end{equation}
where the summation is taken over all spanning trees of $\Gamma$.
A particular consequence is that the computation is independent of the 
choice of root $v$. Note also that, since the degree of $\D_\ell(t)$
coincides with the Euler characteristic of the canonical Seifert
surface, the leading (resp.\ trailing) 
coefficient of $\D_\ell(t)$ is (up to sign) the number of incoherent
(resp.\ coherent) spanning trees in $\Gamma$, and in particular this
number is non-zero. Since $\D_\ell(t)$ is
reciprocal, both quantities coincide. We will for convenience
work with coherent spanning trees. 

\section [Proof of \ref{main}]{Proof of \fullref{main}}

It suffices to consider the
collection ${\cal D}$ of all oriented link diagrams $D$ with twist
number $t(D)$ no greater than an arbitrary integer $n$. Since any
$D\in {\cal D}$ can be obtained from a wiring diagram $X$ of order
$n$, \fullref{th3.4}(1) and (3)  imply that the lengths of $(t+1)^n V_\ell(t)$ and 
$(v^2-1)^{m} (t+1)^{n-m}\hat P_\ell(v,t)$ are bounded, where $m$
is the number of anti-parallel 
twist sites in $X$, and $\ell$ ranges over the links with diagrams
associated to $X$. In the latter case, the set of lengths remains
bounded when we multiply by $(v^2-1)^{m-n}(t+1)^m$. There are finitely
many wiring diagrams of order $n$ with no \tsps, and the first part
of \fullref{main} is proved.

If we perform repeated parallel twisting at one or more sites of an 
oriented link diagram, then the Alexander polynomials of the resulting 
links have bounded Mahler measures by \fullref{lem3.1} together with the previously mentioned fact that $M(p)\le ||p||$, for any nonzero 
polynomial.

However, if we perform repeated anti-parallel
twisting at one or more sites, then, as we will see later in
\fullref{ex5.2}, the euclidean Mahler measures (and hence the Mahler measures) of the Alexander polynomials can grow without bound.
Alternating links have better behavior under such twisting, which we
will explain next.

Consider a wiring diagram of order $n$. Let $D$ be the diagram for the 
link $\ell(q_1, \ldots, q_n)$, the sites $1,\ldots, m$ being
anti-parallel, the remaining sites parallel. 

The form of $(t+1)^{n-m}\D_{\ell(q_1, \ldots, q_n)}(t)$ is given by \fullref{th3.4}(2). As in the case of the Jones polynomial, we can see that the polynomials $$X_{\d_1, \ldots, \d_n}(t^{1/2}, t^{q_{m+1}/2},\ldots, t^{q_n/2})$$ have bounded lengths. To show that the family of 
polynomials $\D_{\ell(q_1, \ldots, q_n)}(t)$ has bounded euclidean Mahler measure, it suffices to show that no coefficient can grow more rapidly as a function of $q_1,\ldots, q_m$ than the leading coefficient. 

Fix $(q_{m+1},\ldots, q_n)\in {\mathbb Z}^{n-m}$ and
$(\d_{m+1},\ldots, \d_n)\in \{0,1\}^{n-m}$, and let $S$
be the set of $m$--tuples ${\d}=(\d_1, \ldots, \d_m)\in \{0,1\}^m$
for which
$$X_\delta(t) = X_{\d_1, \ldots, \d_n}
(t^{1/2}, t^{q_{m+1}/2},\ldots, t^{q_n/2})$$
is non-zero. We will show that if ${\d}$ is maximal in $S$ (in the
sense that ${\d_i'}\ge {\d_i}$, for all $i\le m$, implies that 
${\d'}= {\d}$), then $X_\delta(t)$ and $\D_{\ell(q_1, \ldots, q_n)}(t)$ have the same degree.

First consider the case that $D$ is a special alternating diagram. We make use of the notation established above. Since we are concerned with the growth of coefficients as functions of $q_1, \ldots, q_m$, we may assume that each $q_i$ is at least 2. One easily checks that any anti-parallel twist in $D$ corresponds to a multiple edge in $\Gamma$, while a parallel twist corresponds to a chain (that is, an edge subdivided by any number of vertices). 

It is convenient to form a quotient graph $\bar \Gamma$ by identifying each multiple edge to a single bi-oriented edge, called an {\it anti-parallel edge}.  Edges that are not anti-parallel are said to be {\it ordinary}. We say a tree $\bar T$ in $\bar \Gamma$ to be {\it coherent} if it has a coherent lift in $\Gamma$. 

\begin{figure}[ht!]
\begin{center}
{\includegraphics[width=2.5in]{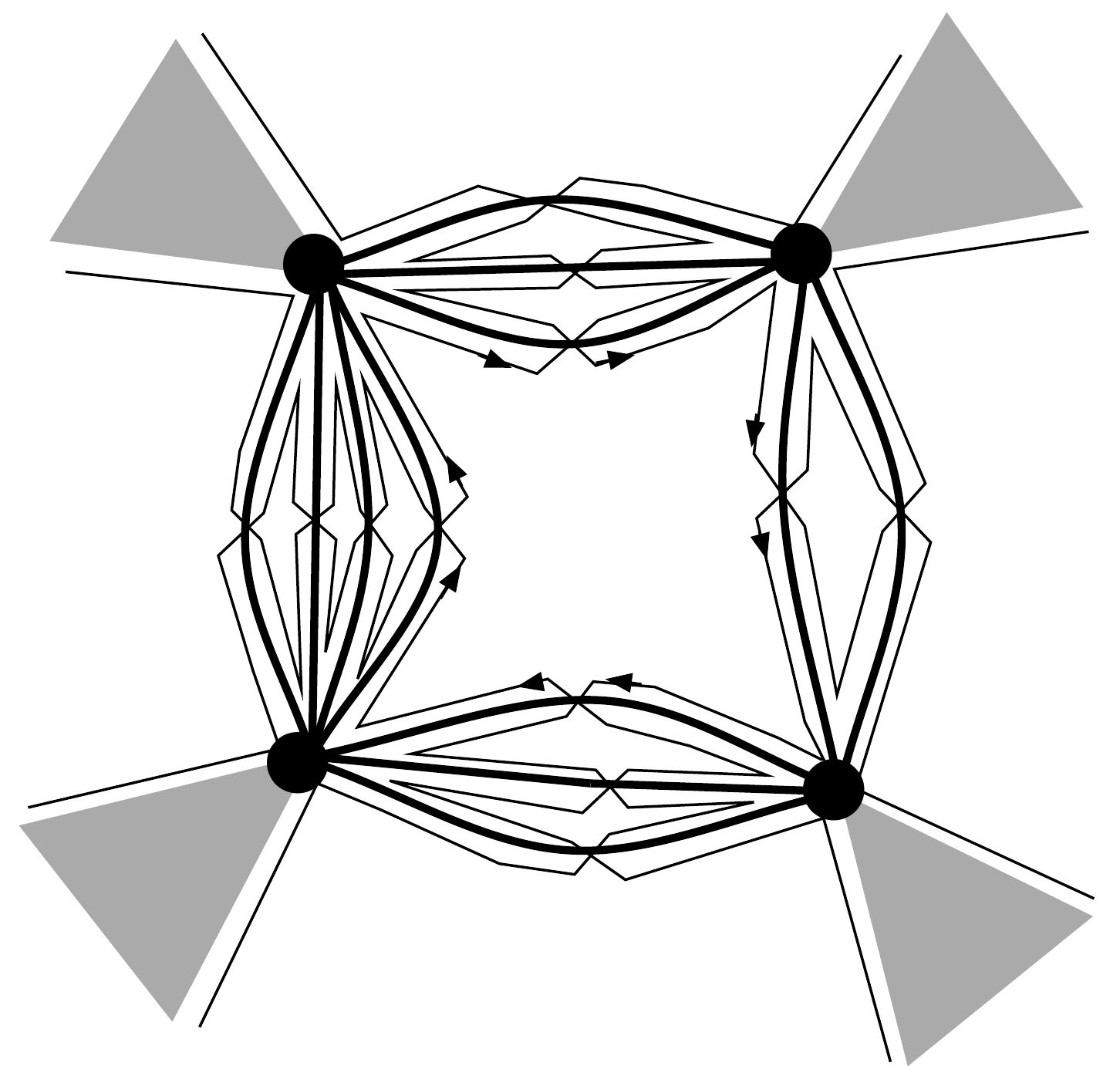}} 
\caption{A cycle of anti-parallel twists} 
\label{twistx}
\end{center}
\end{figure}

Regard $q_1,\ldots, q_m$ as variables, and consider a monomial $q_1^{\d_1}\cdots q_m^{\d_m}$ for which $\d=(\d_1, \ldots, \d_m)$ is maximal. 
The variables $q_i$ that appear (that is, those with exponent $\d_i =1$) correspond to a subgraph $\bar\Sigma$ of anti-parallel edges in $\bar\Gamma$ that contains no cycle, and hence is a forest. For otherwise, the graph $\Gamma$ would have a cycle of multiple edges, as in \fullref{twistx}, and the link $\ell_{(\d_1, \ldots, \d_m)}$, which is obtained by smoothing crossings, would be split.
In that event, $\D_{\ell(\e_1, \ldots, \e_m, q_{m+1}, \ldots, q_n)}$ 
and hence $X_{\d_1, \ldots, \d_m}$ would vanish, contradicting the assumption that the monomial $q_1^{\d_1}\cdots q_m^{\d_m}$ appears. 

By equation \ref{eq3.2} contributions to the leading coefficients of
$\D_{\ell(q_1, \ldots, q_n)}$ from coherent spanning trees do not cancel.
To ensure that $q_1^{\d_1}\cdots q_m^{\d_m}$ appears in this
leading coefficient, we need to show that $\bar\Sigma$ has a
lift $\Sigma$ in $\Gamma$ that is contained in a coherent
spanning tree.

\begin{lem}Any forest $\bar\Sigma \subset \bar \Gamma$ of anti-parallel edges extends to a coherent spanning tree $\bar T$.\end{lem}

\begin{figure}[ht!]\small
\labellist
\pinlabel $v_0$ [bl] at 221 394
\pinlabel $v_1$ [tl] <0pt,2pt> at 225 218
\hair 1.5pt
\pinlabel $e$ [br] at 196 372
\pinlabel $e'$ [l] at 265 312
\hair 8pt
\pinlabel {old tree} [t] at 59 160
\pinlabel {new tree} [t] at 397 160
\pinlabel {ordinary edge} [l] at 164 81
\pinlabel {anti-parallel edge} [l] at 164 34
\endlabellist
\begin{center}
{\includegraphics[width=2.7in]{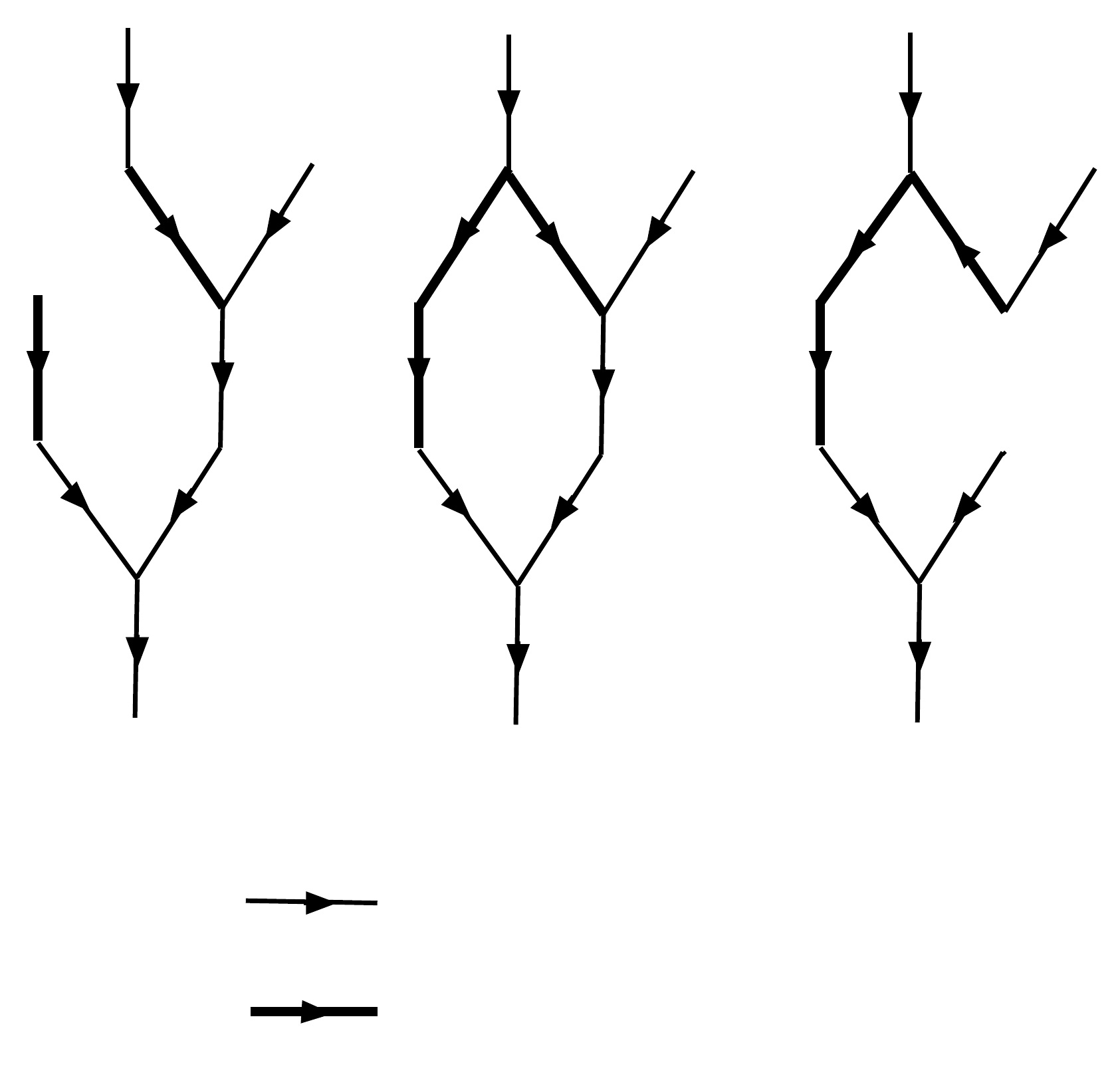}} 
\caption{Exchanging edges $e$ and $e'$ in $\bar T$} 
\label{tree}
\end{center}
\end{figure}

\begin{proof} Let $\bar\Sigma$ be a forest consisting of 
anti-parallel edges. Let $\bar T$ be a coherent spanning tree
(which exists by our preliminary remarks).  By the {\it height} of a
vertex $v$ in $\bar T$ we will mean the distance of that vertex to the root of
$\bar T$.  The root has height $0$.

If some edge $e$ of $\bar\Sigma$ is not in $\bar T$, then we add it to
$\bar T$, thereby creating a unique cycle $\bar C$. We choose an
orientation of $e$ toward the vertex of lower height, if
the heights of the two vertices connecting $e$ differ; otherwise, either orientation for $e$ will do. Consequently, $\bar C$ consists of two oriented paths. In other words, there is a unique vertex $v_0$ in
$\bar C$ with two outgoing and
another vertex $v_1$ with two incoming edges. ($v_0$ is actually
the vertex $e$ points away from.)
Let $x$ be the path in $\bar C$ from $v_0$ to $v_1$
that contains $e$, and let $y$ be the complementary path.

\medskip{\bf Case 1}\qua The path $y$ consists of edges in $\bar \Sigma$.
Following the path $x$ from $v_0$,  let $e'$ be the first edge in $x$ contained in $\bar C\setminus \bar \Sigma.$ Such an edge must exist
since $\bar \Sigma$ is a forest. Delete $e'$ and change the
orientation of  all edges in $x$ that precede $e'$,
including $e$.

\medskip{\bf Case 2}\qua The path $y$ has an edge in $\bar C\setminus \bar \Sigma.$ Let $e'$ be the first such edge in $y$. Delete $e'$ and change the
orientation of all edges in $y$ that precede $e'$. See \fullref{tree}.

\medskip
In this way we obtain a new coherent spanning tree, called $\bar T$ by
abuse of notation, containing one more edge of $\bar \Sigma$ than
before.

Repeat the procedure until $\bar\Sigma$ is
contained in $\bar T$. \end{proof}

General alternating diagrams $D$ can be expressed as  a product $D_1*\cdots *D_r$ of special alternating (Murasugi) factors. As before, we assume that each of $q_1, \ldots, q_m$ is greater than $1$. 

Recall that each anti-parallel twist appears in some factor.  Again consider any monomial  $q_1^{\d_1}\cdots q_m^{\d_m}$ that appears in the leading coefficient of the right-hand side of the equation in \fullref{th3.4}(2). Its variables  can be partitioned so that those in the $i$th subset arise from an anti-parallel twists in the diagram $D_i,\ 1\le i\le r$. Consider such a subset, say $q_1, \ldots, q_{m_1}$ after renumbering. The $m_1$ twists correspond to anti-parallel edges $e_1, \ldots, e_{m_1}$ in the graph $\bar \Gamma$ corresponding to the special factor $D_1$ of $D$. As before, the edges do not form any cycles, and so from them we can form a
coherent spanning tree for $\Gamma$. Consequently, $q_1^{\d_1}\cdots q_m^{\d_{m_1}}$ appears in the leading coefficient of the Alexander polynomial of the Murasugi factor $D_1$, and hence it appears in the leading coefficient of 
$\D_{\ell(q_1, \ldots, q_n)}$. 
To finish the proof we remark again that there are finitely
many wiring diagrams of order $n$ with no \tsps, and
\tsps{} do not alter $M_e(\D)$. This completes the proof of \fullref{main}.

\begin{rem} \fullref{main}(2) is true by the same argument for
the more general class of homogeneous diagrams defined in
\cite{cromwell}. \end{rem}
\begin{rem}
One easily observes that the various bounds on polynomial lengths,
magnitude and number of non-zero coefficients, and Mahler measure
obtained above are (and must be) exponential in the twist number.
We have not elaborated on this theme here, but with the inequalities
in \cite{stoi03} and a bit extra work, it is possible to give explicit
(though likely still not optimal) bounds on the bases of the exponentials. 
\end{rem}

\section{Examples} The set of Mahler measures of Alexander polynomials corresponding to knots obtained from a given knot by twisting repeatedly at several sites may have infinitely many limit points. We illustrate this in \fullref{ex5.1}.

\begin{figure}[ht!]
\begin{center}
\includegraphics[width=2in]{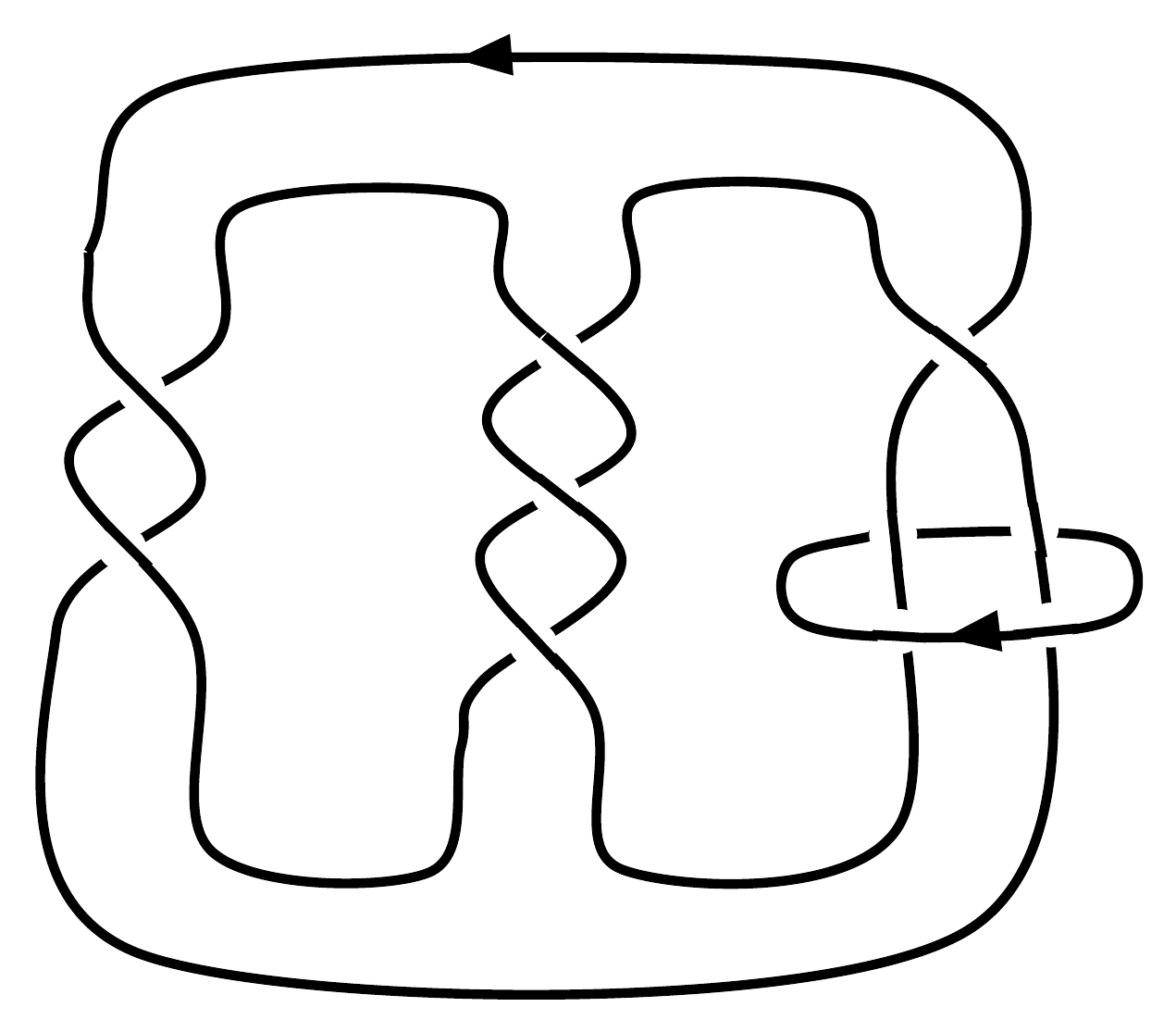} 
\caption{Pretzel link $\ell(2,3, \infty)$} 
\label{pretzel}
\end{center}
\end{figure}

\begin{exa}\label{ex5.1} 
Consider the family of alternating pretzel knots $\ell(2, 2k+1, 2l+1)$
with $2, 2k+1$ and $2l+1$ half-twists in first, second and third
bands, respectively, $k, l\ge 0$.  If we fix $k$ and let $l$ tend
toward infinity, then by Theorem 2.2 of \cite{SW04} the Alexander
polynomials of the resulting knots have Mahler measures that approach
the limit $M(\D_{\ell(2, 2k+1, \infty)}(x,z))$, where $\ell(2, 2k+1,
\infty)$ is the $2$--component link shown in \fullref{pretzel} for the case
$k=1$.

Using Remark 3.3 of \cite{SW04}, we can compute $\D_{\ell(2, 2k+1, \infty)}(x,z)$
from the 3--variable Alexander polynomial $\D_{\ell(2, \infty, \infty)}(x,y,z)$, where  $\ell(2, \infty, \infty)$ is obtained from\break $\ell(2, 2k+1, \infty)$ by leaving a single half-twist in the second band and encircling it, just as we have done with the third band. We have
$$(x^2-1) \D_{\ell(2, 2k+1, \infty)}(x,z)\,=\,
 \D_{\ell(2, \infty, \infty)}(x,x^{-2k},z)\,.$$
Since the two polynomials differ only by a cyclotomic factor,
$$M(\D_{\ell(2, 2k+1, \infty)}(x,z))=
M(\D_{\ell(2, \infty, \infty)}(x,x^{-2k},z)).$$
The polynomial $\D_{\ell(2, \infty, \infty)}
(x,x^{-2k},z)$ can be expressed, up to multiplication by a unit, as
$$(x-1)\ \Big [\ (1-x-x^2-x^{2k+2}) +
x^{2k+4}z(1-x^{-1}-x^{-2}-x^{-2k-2})\ \Big]. $$ A technique of D Boyd
\cite{boyd} (see \cite[Lemma 4.1]{SW04}) enables us to compute the
Mahler measure of this polynomial as the Mahler measure of the
single-variable polynomial $f_k(x)=1-x-x^2-x^{2k+2}$. An argument
based on Rouch\'e's theorem (Brown and Churchill \cite{BC}) shows
that $f_k(x)$ has exactly one zero $\zeta_k$ inside the unit circle, a
zero that is real.  The product of all the moduli of zeros of $f_k$ is
equal to the modulus of the absolute term $1$, so
$M(f_k)=1/|\zeta_{k}|$.  Since for $k\ne k'$ the difference
$f_k(x)-f_{k'}(x)$ is a unit times a cyclotomic polynomial, $f_k$ and
$f_{k'}$ have no common zeros off the unit circle. In particular,
$\zeta_k \ne \zeta_{k'}$ whenever $k\ne k'$. Now at most two real
numbers have the same modulus, thus the Mahler measures of the
polynomials $f_k$ are triplewise distinct, and so the set of values
$M(\D_{\ell(2, 2k+1, 2l+1)}(x))$ has infinitely many distinct limit
points.

\fullref{main} implies that the limit points are bounded. \end{exa}

\fullref{ex5.2} shows that the conclusion of \fullref{main}(2) does not hold if the hypothesis that the diagrams be alternating is dropped. 

\begin{figure}[ht!]\small
\labellist
\pinlabel $n$ [r] at 67 297
\endlabellist
\begin{center}
{\includegraphics[width=2in]{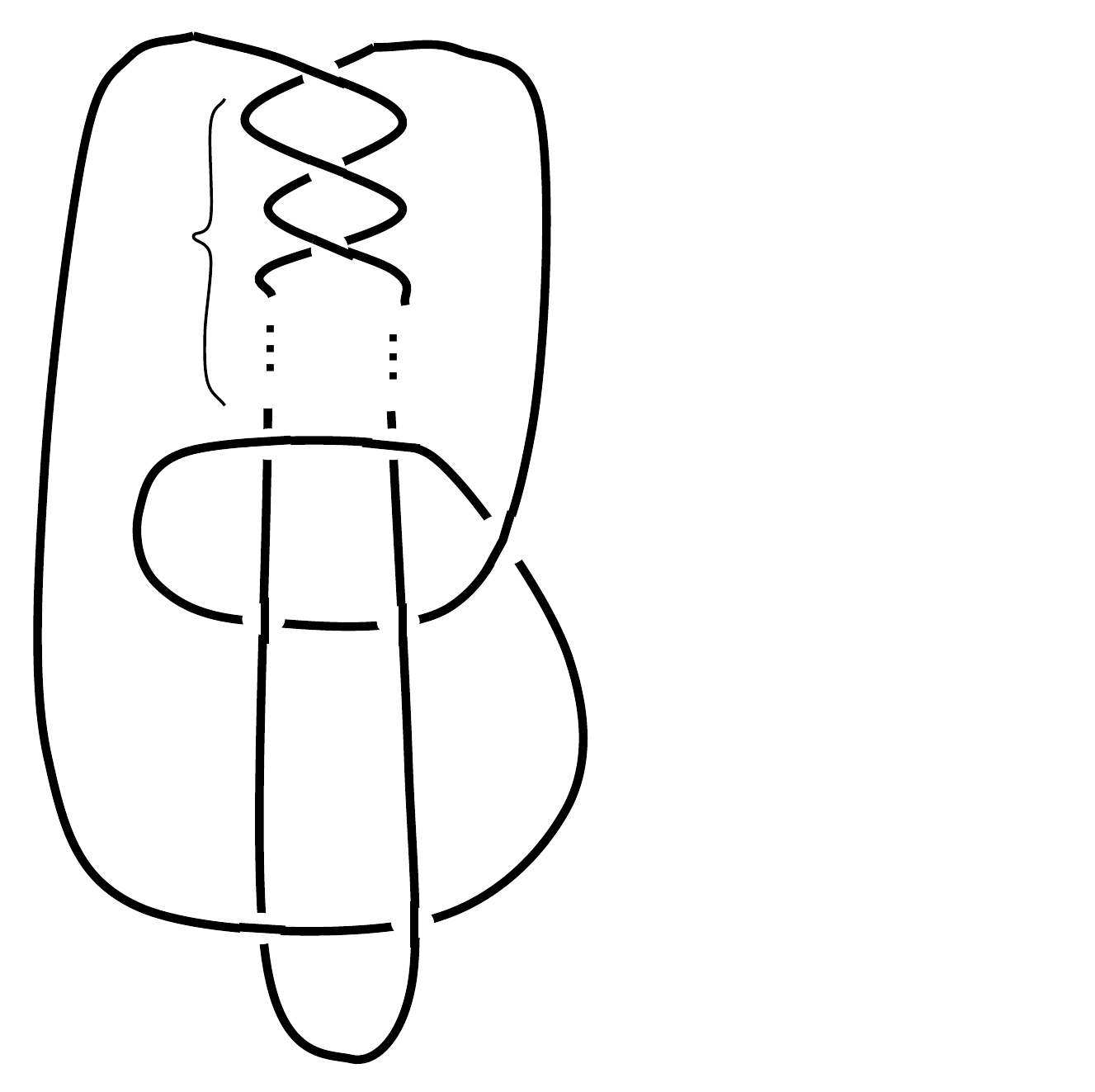}} 
\caption{Non-alternating knot diagram $D_n$} 
\label{twist3}
\end{center}
\end{figure}

\begin{exa}\label{ex5.2} Consider the
family of non-alternating diagrams $D_n$ in \fullref{twist3}, each member containing $n=1,3,5,\ldots$ half-twists in its upper portion. 
The Alexander polynomials
of the corresponding knots $k_n$ are easily seen to be
$$\D_{k_n}(t)= t^4- {n+5\over 2} t^3 + (n+4) t^2 - {n+5\over 2} t +1.$$
The twist numbers of the diagrams $D_n$ are all equal. However, the euclidean Mahler measures of the polynomials must tend to infinity, since the sum of the four roots is $(n+5)/2$ (and so one root has modulus at
least $(n+5)/8$). Since the volumes of $k_n$ are also bounded (by
Thurston's hyperbolic surgery theorem; see Lackenby \cite{lacken}), similarly
\fullref{cor2.3} also fails without the alternation assumption. \end{exa}

We conclude by showing that the converse of \fullref{main}(2) does not hold. 

\begin{exa}Consider the tangle $T$ in \fullref{twist4}. By replacing a single crossing  by $T$ in the diagram of a trefoil, as in the figure, we obtain the knot $8_{10}$. Iterating the procedure, always replacing  (for the sake of definiteness) the top crossing, we obtain a sequence of diagrams $D_n$ for alternating knots $k_1 = 3_1, k_2 = 8_{10}$, et cetera.

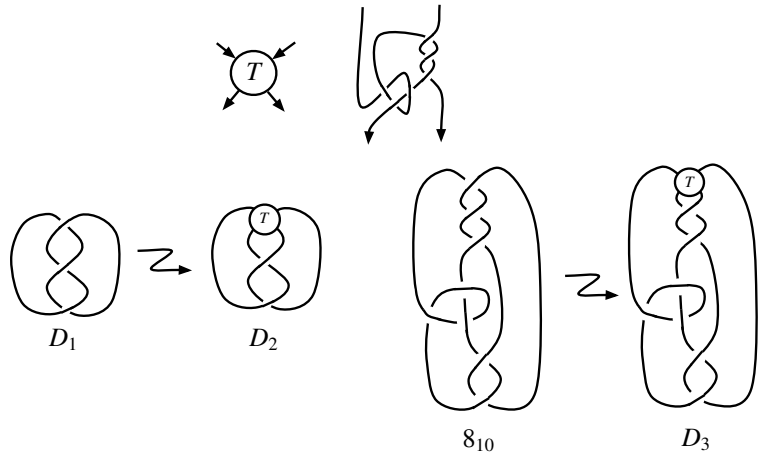
\begin{figure}[ht!]\small
\labellist
\pinlabel $T$ <1pt,1pt> at 320 566
\pinlabel {$\scriptscriptstyle T$} <-.2pt,1pt> at 337 379
\pinlabel {$\scriptscriptstyle T$} <1pt,.5pt> at 876 427
\pinlabel $D_1$ at 80 227
\pinlabel $D_2$ at 336 227
\pinlabel $8_{10}$ at 610 99
\pinlabel $D_3$ at 888 99
\endlabellist
\begin{center}
{\includegraphics[width=4in]{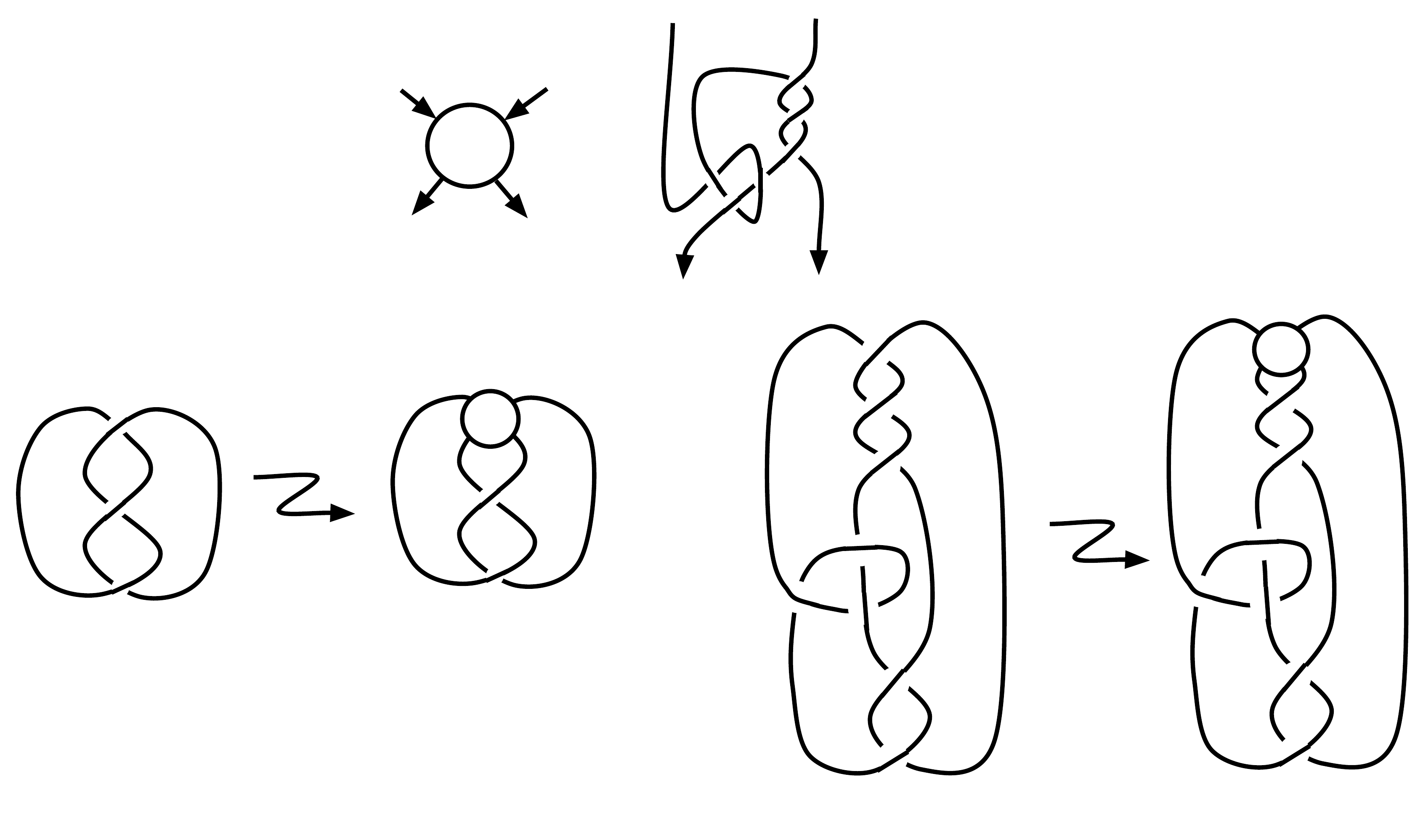}} 
\caption{Tangle $T$ and alternating knot diagrams} 
\label{twist4}
\end{center}
\end{figure}

We claim that $\D_{k_n}(t) = (t-1+t^{-1})^{2n-1}$, for each $n$. Consider the Conway skein module over the field of fractions of
${\mathbb Z}[t^{1/2}, t^{-1/2}]$ generated by oriented tangles modulo 
the relations
$\ell_+ - \ell_- -(t^{1/2}-t^{-1/2})\ell_0$, where $(\ell_+, \ell_-, \ell_0)$ is any skein triple. 
The tangle $T$ induces a mapping $F:S \mapsto F(S)= \D_{D(T\cdot S)}(t)$, where $D(T\cdot S)$ denotes the denominator closure of the tangle product of $T$ and $S$ (see \fullref{twist5}). 

\begin{figure}[ht!]\small
\labellist
\pinlabel $T$ <0pt,1pt> at 305 218
\pinlabel $S$ <0pt,1pt> at 305 114
\pinlabel $S$ <0pt,1pt> at 57 165
\pinlabel $F(S)$ [l] at 487 147
\pinlabel {$D(T\cdot S)$} <0pt, -2pt> at 305 40
\endlabellist
\begin{center}
{\includegraphics[width=3in]{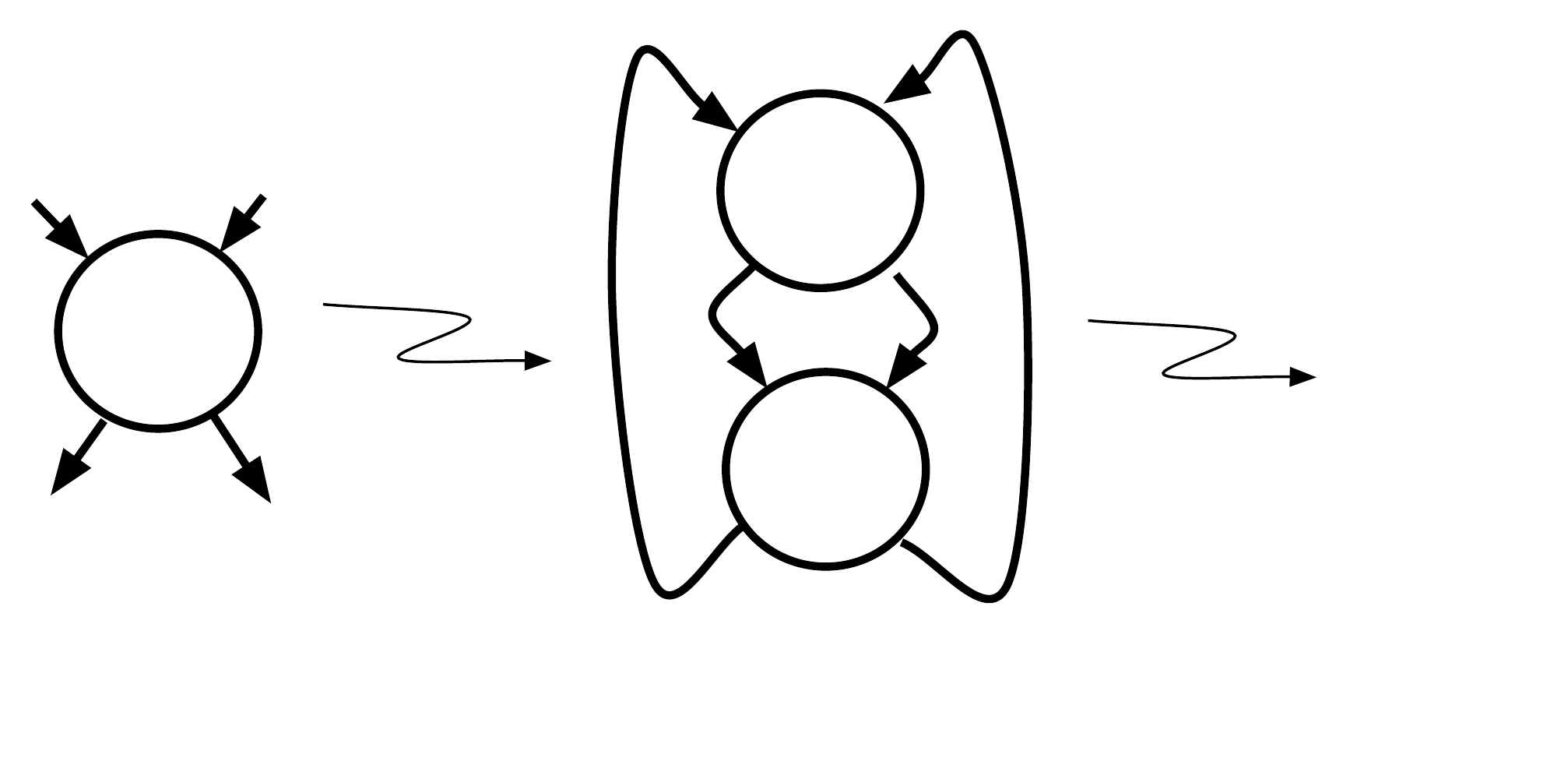}} 
\caption{Linear functional $F$} 
\label{twist5}
\end{center}
\end{figure}
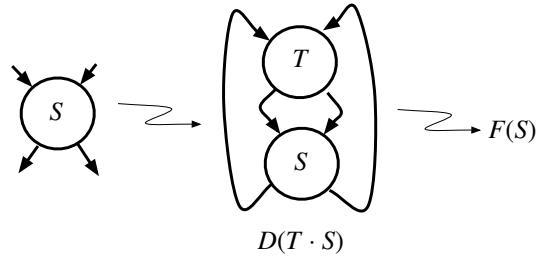

The tangles $S_1, S_2$ in \fullref{twist7} form a basis for the skein module. We can write $T = f(t)S_1+ g(t)S_2,$ for some scalars $f(t),g(t)\in {\mathbb Z}[t^{1/2}, t^{-1/2}]$.  Then
\begin{equation}
F(S) = f(t)\cdot \D_{D(S_1\cdot S)}(t) + g(t)\cdot \D_{D(S_2\cdot S)}(t),
\label{eq5.1}
\end{equation}
for any tangle $S$.

\begin{figure}[ht!]\small
\labellist\hair 10pt
\pinlabel $S_1$ [t] at 35 72
\pinlabel $S_2$ [t] at 200 72
\endlabellist
\begin{center}
{\includegraphics[width=1.5in]{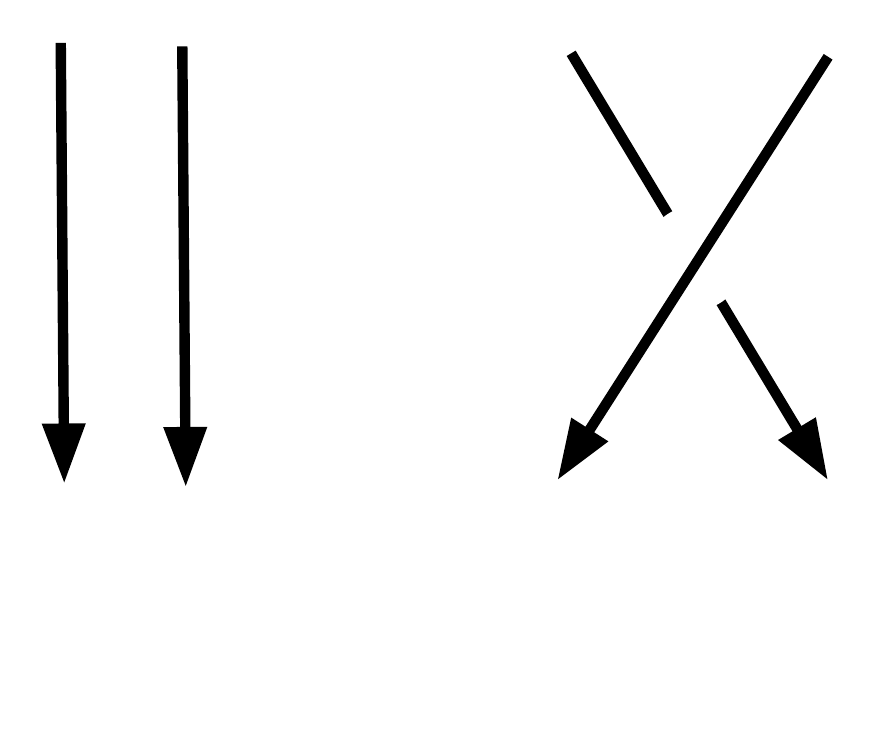}} 
\caption{ Basis $S_1, S_2$} 
\label{twist7}
\end{center}
\end{figure}

We can find $f(t)$ and $g(t)$ easily by substituting proper values
for $S$ in \ref{eq5.1}. If $S=S_1$, then $D(T\cdot S)$ is the square knot,
$D(S_1\cdot S)$ is a trivial 2--component link, and $D(S_2\cdot S)$ is an
unknot. Hence 
$(t-1+t^{-1})^2 = f(t)\cdot 0 + g(t)\cdot 1$, and so $g(t) = (t-1+t^{-1})^2$.
Similarly, 
if $S= S_2\cdot S_2$, then $D(T\cdot S)$ is the knot $8_{10}$, $D(S_1\cdot S)$ is a Hopf link, and $D(S_2\cdot S)$ is a trefoil. So
\begin{equation*}\begin{split}
\D_{8_{10}}& =(t-1+t^{-1})^3\,=\,(t^{1/2}-t^{-1/2})\cdot f(t)\,+\,
(t-1+t^{-1})\cdot g(t)  \cr
&=\,(t^{1/2}-t^{-1/2})\cdot f(t)\,+\,(t-1+t^{-1})^3\,.
\end{split}\end{equation*}
Thus $f(t)=0$, and in the skein module 
$T=(t-1+t^{-1})^2 S_2$, and then $\D_{k_n}=(t-1+t^{-1})^2
\D_{k_{n-1}}$.  Since $\D_{k_1}(t)= t-1+t^{-1}$, we have $\D_{k_n}(t)= (t-1+t^{-1})^{2n-1}$. 

Observe that the diagrams in \fullref{twist4} have unbounded twist number (and
hence, by \cite{lacken} also volume), while the Alexander polynomials
$\D_{k_n}(t)$ are all products of cyclotomic polynomials and hence
have trivial Mahler measure. \end{exa}

\begin{rem} If one drops the alternation assumption, then
for any knot $k$ there exists a hyperbolic knot $\tilde k$ with
arbitrarily large volume (and twist number) and the {\sl same}\/ 
Alexander polynomial as $k$. For trivial polynomial knots this result
was proved in \cite{kalfa}.  Another method that constructs for any
(admissible) Alexander knot polynomial (in fact, Alexander
invariant) an arborescent knot $\tilde k$ is given by the Stoimenow
in \cite{stoi05}. Silver and Whitten \cite{SWh} have used a
yet different construction to choose $\tilde k$ so that the group of
$\tilde k$ can be mapped onto the group of $k$ sending meridian to
meridian and longitude to longitude. \end{rem}

\bibliographystyle{gtart}
\bibliography{link}

\end{document}